\documentclass[12pt,leqno]{amsart}
\newtheorem{theo}{Theorem}[section]

\newtheorem{pr}{Proposition}[section]
\newtheorem{rem}{Remark}[section]
\newcommand{\be}{\begin{equation}}
\newcommand{\ee}{\end{equation}}
\newcommand{\bea}{\begin{eqnarray}}
\newcommand{\eea}{\end{eqnarray}}
\newcommand{\beb}{\begin{eqnarray*}}
\newcommand{\eeb}{\end{eqnarray*}}
\usepackage{amssymb,amsfonts,amsthm,setspace,indentfirst,mathrsfs}
\numberwithin{equation}{section}
\usepackage[dvips]{graphics}
\usepackage{epsfig}

\numberwithin{figure}{section}
\begin{document}
\title[On the existence of a generalized class of recurrent manifolds]{ {\normalsize \bf{On the existence of a generalized class of recurrent manifolds}}}
\author[A. A. Shaikh, I. Roy and H. Kundu]{Absos Ali Shaikh, Indranil Roy and Haradhan Kundu}
\date{\today}
\address{\noindent\newline Department of Mathematics,\newline University of
Burdwan, Golapbag,\newline Burdwan-713104,\newline West Bengal, India}
\email{aask2003@yahoo.co.in, aashaikh@math.buruniv.ac.in}
\email{royindranil1@gmail.com}
\email{kundu.haradhan@gmail.com}
\dedicatory{Dedicated to Prof. Dr. Constantin Udri\c{s}te on his 76th birthday}
\begin{abstract}
The present paper deals with the proper existence of a generalized class of recurrent manifolds, namely, hyper-generalized recurrent manifolds. We have established the proper existence of various generalized notions of recurrent manifolds. For this purpose we have presented a metric and computed its curvature properties and finally we have obtained the existence of a new class of semi-Riemannian manifolds which are non-recurrent but hyper-generalized recurrent, Ricci recurrent, conharmonically recurrent, manifold of recurrent curvature 2-forms and semisymmetric; not weakly symmetric but weakly Ricci symmetric, conformally weakly symmetric and conharmonically weakly symmetric; non-Einstein but Ricci simple; do not satisfy $P\cdot P =0$ but fulfills the condition $P \cdot P = -\frac{1}{3}Q(S, P)$, $P$ being the projective curvature tensor.
\end{abstract}
\maketitle
\flushbottom
\section{\bf Introduction}
Let $M$ (dim $M = n \ge 3$) be a connected semi-Riemannian manifold endowed with a semi-Riemannian metric $g$ of signature
 ($s, n-s$), $0\le s \le n$. If $s=0$ or $n$ then $M$ is a Riemannian manifold, and if $s=1$ or $n-1$ then $M$ is a Lorentzian
 manifold. A semi-Riemannian manifold has mainly three notions of curvature, namely, the Riemann-Christoffel curvature tensor $R$ (simply called curvature tensor), the Ricci tensor $S$ and the scalar curvature $r$. Let $\nabla$ be the Levi-Civita connection on $M$, which is the unique torsion free metric connection on $M$.\\
\indent It is well known that $M$ is locally symmetric if its local geodesic symmetries at each point are all isometry. In terms of curvature restriction, $M$ is locally symmetric if $\nabla R = 0$ (\cite{Ca26}, \cite{Ca27}). The study of generalization of locally symmetric manifolds was started in 1946 and continued till date in different directions, such as semisymmetric manifolds by Cartan \cite{Ca46} (which were latter classified by Szab$\acute{\mbox{o}}$ (\cite{Sz82, Sz84, Sz85})), $\kappa$-spaces or recurrent spaces by Ruse (\cite{Ru46}, \cite{Ru49}, \cite{Ru49a}), weakly symmetric manifolds by Selberg \cite{Se56}, pseudosymmetric manifolds by Chaki \cite{Ch87}, pseudosymmetric manifolds by Deszcz (\cite{Des92}, \cite{DGr87}), weakly symmetric manifolds by Tam$\acute{\mbox{a}}$ssy and Binh \cite{TB89} etc.\\
\indent Again in the literature of differential geometry we find various curvature restricted geometric structures formed by imposing a restriction on the Ricci tensor, such as Ricci symmetric manifold ($\nabla S = 0$), Ricci recurrent manifold by Patterson \cite{Pa52}, Ricci semisymmetric manifold by Szab$\acute{\mbox{o}}$ (\cite{Sz82, Sz84, Sz85}), Ricci pseudosymmetric manifold by Deszcz \cite{Des89}, pseudo Ricci symmetric manifold by Chaki \cite{Ch88}, weakly Ricci symmetric manifold by Tam$\acute{\mbox{a}}$ssy and Binh \cite{TB93} etc.\\
\indent During the study of simply harmonic Riemannian spaces, Ruse (\cite{Ru46}, \cite{Ru49}, \cite{Ru49a}) introduced the notion of kappa spaces, which were latter known as recurrent spaces (see, \cite{Ru49a}, \cite{Wa50}). The manifold $M$ is said to be recurrent (\cite{Ru46}, \cite{Ru49}, \cite{Ru49a}) (resp., Ricci recurrent \cite{Pa52}) if
\[\nabla R = A \otimes R \ \ \ \left(\mbox{resp.} \ \nabla S = A \otimes S \ \right)\]
holds for some 1-form $A$. During the last six decades, in literature we find various generalized notions of recurrent manifolds, such as 2-recurrent manifolds by Lichnerowicz \cite{Li52}, conformally recurrent manifolds by Adati and Miyazawa \cite{AM67}, projectively recurrent manifolds by Adati and Miyazawa \cite{AM77},  concircular recurrent or generalized recurrent manifolds by Dubey \cite{Du79}, quasi-generalized recurrent manifolds by Shaikh and Roy \cite{SR10}, hyper-generalized recurrent manifolds by Shaikh and Patra \cite{SP10}, weakly generalized recurrent manifolds by Shaikh and Roy \cite{SR11}, super generalized recurrent manifolds by Shaikh et al. \cite{SRK16},  manifolds of recurrent curvature 2-forms (\cite{Be87}, \cite{LR89}) etc. Recently, Olszak and Olszak \cite{OO12} showed that every generalized recurrent manifold is recurrent (see also \cite{MIKES76}, \cite{MIKES96}, \cite{MVH09}). It is interesting to note that although the notion of generalized recurrent manifold does not exist but the proper existence of weakly generalized recurrent (resp., quasi-generalized recurrent and super generalized recurrent) manifold is presented in \cite{SAR13} (resp., \cite{SRK16}).\\
\indent The study of recurrence condition is a study of tensors which are invariant (upto a scaling) under a parallel transport. It is an important topic of research as the recurrent condition implies that the curvature tensor remains invariant under the covariant differentiation \cite{Ru49a}. Recent investigation of Hall and Kirik \cite{HK15} shows that the study of recurrent manifolds is closely related to the holonomy theory.\\ 
\indent As a generalization of recurrent manifolds, recently Shaikh and Patra \cite{SP10} introduced the notion of hyper-generalized recurrent manifolds and studied its various geometric properties. The present paper is devoted to the study of the proper existence of hyper-generalized recurrent manifolds and the interrelation of such notion with other geometric structures. The paper is organized as follows. Section 2 deals with various curvature restricted geometric structures and their interrelations as preliminaries. The last section is devoted to the proper existence of hyper-generalized recurrent manifold. For this purpose we pose the following question:\\
\noindent\textbf{Q. 1} {\it Does there exist a hyper-generalized recurrent manifold which is\\
(i) not locally symmetric but semisymmetric,\\
(ii) not weakly symmetric but weakly Ricci symmetric as well as conformally weakly symmetric,\\ 
(iii) not recurrent but Ricci recurrent, conformally recurrent as well as a manifold of recurrent curvature 2-forms,\\
(iv) neither class $\mathcal{A}$ nor $\mathcal{B}$ but of constant scalar curvature,\\
(v) not Einstein but Ricci simple i.e., special quasi-Einstein,\\
(vi) does not satisfy $P\cdot P =0$ but fulfills the condition $P\cdot P = L Q(S, P)$ for some scalar $L$.}
\vspace{0.1in}\\
\indent In the last section of this paper we have computed the curvature properties of the metric given by
\[ds^2 = g_{ij} dx^i dx^j = e^{x^1+x^3} (dx^1)^2 + 2 dx^1 dx^2 + (dx^3)^2 + e^{x^1} (dx^4)^2, \ \ i, j = 1, 2, 3, 4,\]
which produces the answer of \textbf{Q. 1} as affirmative, and finally we obtain the existence of a new class of semi-Riemannian manifolds (see Theorem \ref{th3.1} and Remark \ref{rem3.1}).
%
\section{Preliminaries}
\indent Let us consider a connected semi-Riemannian manifold $M$ of dimension $n (\ge 3)$ (these conditions are assumed throughout the paper). The manifold $M$ is said to be Einstein (briefly, $E_n$) if its Ricci tensor $S$ of type $(0, 2)$ satisfies $S =\frac{r}{n} g$. The manifold is said to be quasi-Einstein (briefly, $QE_n$) (\cite{DDHKS00}, \cite{DDVV99}, \cite{DGHS98}, \cite{DHS01}, \cite{DHS01a}, \cite{DVV96}, \cite{Gl08}, \cite{SH11}, \cite{SKH11}, \cite{SYH09}) if
\[S - \alpha g = \beta \eta \otimes \eta\]
 on $\{ x \in M: \left(S - \frac{r}{n}\, g\right)_x \neq 0\}$, where $\otimes$ is the tensor product, $\eta$ is an 1-form and $\alpha, \beta$ are two scalars on that set. It may be noted that if $\alpha$ vanishes identically then a quasi-Einstein manifold turns into a Ricci simple manifold. Hence a semi-Riemannian manifold $(M, g)$, $n \ge 3$, is said to be Ricci simple if the following condition
\be\label{2.1}
S = \beta \eta \otimes \eta
\ee
holds for a scalar $\beta$ and a 1-form $\eta$.
\begin{pr}\label{pr2.1}
A manifold satisfying \eqref{2.1} is Ricci recurrent if $\eta$ is parallel (i.e., $\nabla \eta =0$).
\end{pr}
\begin{pr}\label{pr2.2}
A manifold satisfying $R = \alpha~ g\wedge(\eta\otimes \eta)$ for a scalar $\alpha$ and a 1-from $\eta$, then the manifold is recurrent if $\eta$ is parallel. The result is also true for any other (0,4)-tensor.
\end{pr}
\indent The manifold $M$ is said to be of class $\mathcal B$ or of Codazzi type Ricci tensor (\cite{Fe81}, \cite{Si81}) (resp. class $\mathcal A$ or cyclic Ricci parallel \cite{Gr78}) if
\[(\nabla_{X_1} S)(X_2, X_3) = (\nabla_{X_2} S)(X_1, X_3)\]
\[\left(\mbox{resp.} \ \ \ (\nabla_{X_1} S)(X_2, X_3) + (\nabla_{X_2} S)(X_3, X_1) + (\nabla_{X_3} S)(X_1, X_2) = 0 \ \right)\]
holds for all vector fields $X_1, X_2 , X_3 \in \chi(M)$, where $\chi(M)$ being the Lie algebra of all smooth vector fields on $M$. Throughout this paper we consider $X, Y, X_i\in \chi(M)$, $i=1,2,3,\cdots$.\\
\indent It is well known that the conformal transformation is an angle preserving mapping, the projective transformation is a geodesic preserving mapping whereas concircular transformation is the geodesic circle preserving mapping and conharmonic transformation is the harmonic function preserving mapping and each transformation induces a curvature tensor as invariant, namely, conformal curvature tensor $C$, projective curvature tensor $P$, concircular curvature tensor $W$ and conharmonic curvature tensor $K$ (\cite{DDHKS00}, \cite{DGHS11}, \cite{Gl08}, \cite{Is57}, \cite{YK89}), and are respectively given as:
\beb
C(X_1,X_2,X_3,X_4)= \left[R-\frac{1}{n-2}(g\wedge S) + \frac{r}{2(n-2)(n-1)}(g\wedge g)\right](X_1,X_2,X_3,X_4),
\eeb
\beb
P(X_1, X_2, X_3, X_4) &=& R(X_1, X_2, X_3, X_4)\\
&& - \frac{1}{n-1}[g(X_1, X_4)S(X_2, X_3)-g(X_2, X_4)S(X_1, X_3)],\\
W(X_1,X_2,X_3,X_4)&=& \left[R-\frac{r}{2n(n-1)}(g\wedge g)\right](X_1,X_2,X_3,X_4),\\
K(X_1,X_2,X_3,X_4)&=& \left[R-\frac{1}{n-2}(g\wedge S)\right](X_1,X_2,X_3,X_4),
\eeb
where `$\wedge$' denotes the the Kulkarni-Nomizu product. For two $(0,2)$-tensors $J$ and $F$ the Kulkarni-Nomizu product is defined as (see e.g. \cite{DGHS11}, \cite{Gl08}, \cite{Ko06}):
\begin{eqnarray*}
(J\wedge F)(X_1,X_2,X_3,X_4) &=& J(X_1,X_4)F(X_2,X_3) + J(X_2,X_3)F(X_1,X_4)\\
&-& J(X_1,X_3)F(X_2,X_4) - J(X_2,X_4)F(X_1,X_3).
\end{eqnarray*}
\indent We recall that the manifold $M$ is locally symmetric (resp., Ricci symmetric) (\cite{Ca26}, \cite{Ca27}, \cite{Ca46}) if $\nabla R =0$ (resp., $\nabla S =0$). Similarly for other curvature tensors we have conformally symmetric, projectively symmetric, concircularly symmetric, conharmonically symmetric manifolds. We note that local symmetry, projective symmetry and concircular symmetry are equivalent (see \cite{SK14} and also references therein). An $n$-dimensional locally symmetric (resp., Ricci symmetric and conformally symmetric) manifold is denoted by $S_n$ (resp., $RS_n$ and $CS_n$).\\
\indent Again the manifold $M$ is said to be recurrent (\cite{Ru46}, \cite{Ru49}, \cite{Ru49a}, \cite{Wa50}) (resp., Ricci recurrent \cite{Pa52}) if 
\[\nabla R = A \otimes R \ \ \ \left(\mbox{resp.} \ \nabla S = A \otimes S \ \right)\]
holds for a 1-form $A$ on $\{x\in M:R_x\neq 0\}$ (resp., $\{x\in M: S_x\neq 0\}$). Similarly we can define the conformally recurrent, projectively recurrent, concircularly recurrent, conharmonically recurrent manifolds. Recently Shaikh and Kundu \cite{SK14} showed that the recurrent, projectively recurrent and concircularly recurrent manifolds are equivalent. An $n$-dimensional recurrent (resp., Ricci recurrent and conformally recurrent) manifold is denoted by $K_n$ (resp., $RK_n$ and $CK_n$).\\
\indent The manifold $M$ is said to be quasi-generalized recurrent \cite{SR10} (briefly, $QGK_n$), hyper generalized recurrent \cite{SP10} (briefly, $HGK_n$), weakly generalized recurrent \cite{SR11} (briefly, $WGK_n$) and super generalized recurrent \cite{SRK16} (briefly, $SGK_n$) respectively if
\[\nabla R = A \otimes R + B \otimes g \wedge[g + (\eta \otimes \eta)],\]
\[\nabla R = A \otimes R + B \otimes (S \wedge g),\]
\[\nabla R = A \otimes R + B \otimes \frac{1}{2}(S \wedge S) \mbox{ and}\]
\[\nabla R = A \otimes R + B \otimes (S \wedge S) + D \otimes (g \wedge S) + E \otimes (S \wedge S)\]
holds on 
 $\{x\in M: R_x \neq 0 \mbox{ and }  [g\wedge (g + \eta\otimes\eta)]_x \neq 0\}$,
 $\{x\in M: R_x \neq 0 \mbox{ and }  (g\wedge S)_x \neq 0\}$,
 $\{x\in M: R_x \neq 0 \mbox{ and }  (S\wedge S)_x \neq 0\}$ and
 $\{x\in M: R_x \neq 0 \mbox{ and either }  (S\wedge S)_x  \neq 0 \mbox{ or }  (g\wedge S)_x \neq 0\}$ respectively for some 1-forms $A, B, D, E$ and $\eta$.\\
\indent We note that on a quasi-Einstein manifold with associated scalars $\alpha$ and $\beta$, if $\alpha = \beta$ then the notion of $WGK_n$ and $QGK_n$ are equivalent, and if $2\alpha = \beta$ then $HGK_n$ and $QGK_n$ are equivalent.\\
\indent Again the manifold $M$ is said to be generalized Ricci recurrent if $\nabla S = A \otimes S + B \otimes g$ holds on $\{x\in M: S_x \neq 0\}$ for some 1-forms $A$ and $B$.\\
\indent During the study of a certain kind of conformally flat Riemannian space of class one, Sen and Chaki \cite{SC67} obtained a curvature condition which was also studied by Chaki \cite{Ch87} and named it as pseudosymmetric manifold. The manifold $M$ is said to be pseudosymmetric in sense of Chaki or simply Chaki pseudosymmetric \cite{Ch87} (briefly, $PS_n$) if 
\beb
&&(\nabla_{X} R)(X_1, X_2, X_3, X_4) = 2A(X) R(X_1, X_2, X_3, X_4) + A(X_1) R(X_1, X, X_3, X_4)\\
&& \hspace{2cm} + A(X_2) R(X_1, X, X_3, X_4) + A(X_3) R(X_1, X_2, X, X_4) + A(X_4) R(X_1, X_2, X_3, X)
\eeb
holds $\forall~ X, X_i \in \chi(M)$ and a non-zero 1-form $A$ on $\{x\in M: R_x \neq 0\}$. Similarly we can define the conformally (briefly, $CPS_n$) (resp., projectively, concircularly, conharmonically) pseudosymmetric manifolds. We note that every recurrent manifold is Chaki pseudosymmetric \cite{Ewer93} but the result is not true for other curvature tensors (e.g., see Theorem \ref{th3.1} for the existence of $CK_n$ which is not $CPS_n$).\\ 
\indent Again in 1988, Chaki \cite{Ch88} introduced the notion of pseudo Ricci symmetric manifolds. The manifold $M$ is said to be pseudo Ricci symmetric \cite{Ch88} (briefly, $PRS_n$) if
\[
(\nabla_X S)(X_1,X_2) = 2A(X) S(X_1,X_2) + A(X_1) S(X, X_2) + A(X_2) S(X_1,X)
\]
holds $\forall~ X, X_i \in \chi(M)$ and a non-zero 1-form $A$ on $\{x\in M: S_x \neq 0\}$.\\
\indent As a generalization of Chaki pseudosymmetric manifold and recurrent manifold, in 1989 Tam$\acute{\mbox{a}}$ssy and Binh \cite{TB89} introduced the notion of weakly symmetric manifold. The manifold $M$ is said to be weakly symmetric by Tam$\acute{\mbox{a}}$ssy and Binh \cite{TB89} (briefly, $WS_n$) if
\beb
&&\nabla_X  R(X_1, X_2, X_3, X_4) = A(X) R(X_1, X_2, X_3, X_4) + B(X_1) R(X_1, X_2, X_3, X_4)\\
&& \hspace{2cm} + \overline B(X_2) R(X_1, X, X_3, X_4) + D(X_3) R(X_1, X_2, X, X_4) + \overline D(X_4) R(X_1, X_2, X_3, X)
\eeb
holds $\forall~ X, X_i \in \chi(M)$ $(i =1,2,3,4)$ and some 1-forms $A, B, \overline B, D$ and $\overline D$  on $\{x\in M: R_x \neq 0\}$. Similarly for other $(0,4)$-tensors we obtain various weak symmetry structures, e.g., conformally weak symmetry (briefly, $CWS_n$), projectively weak symmetry, concircularly weak symmetry, conharmonically weak symmetry etc. It may be mentioned that every weakly symmetric manifold is Chaki pseudosymmetric (\cite{Ewer93}, \cite{SDHJK15}). For the reduced defining condition of various weak symmetry we refer the reader to see the recent investigation on pseudosymmetric manifolds by Shaikh et al. \cite{SDHJK15} and also references therein.\\
\indent Again in 1993 Tam$\acute{\mbox{a}}$ssy and Binh \cite{TB93} introduced the notion of weakly Ricci symmetric manifold. The manifold $M$ is said to be weakly Ricci symmetric \cite{TB93} (briefly, $WRS_n$) if 
\[
(\nabla_{X} S)(X_1,X_2) = A(X) S(X_1,X_2) + B(X_1) S(X, X_2) + D(X_2) S(X_1, X)
\]
holds $\forall~ X, X_i \in \chi(M)$ and some 1-forms $A, B, D$ on $\{x\in M: S_x \neq 0\}$. Again replacing $S$ by some other $(0,2)$-tensors, we get various weak symmetry structures. In \cite{SDHJK15} Shaikh et al. presented the reduced form of the defining condition of weak symmetry for various $(0,2)$-tensors. Weakly symmetric and weakly Ricci symmetric manifolds by Tam\'assy and Binh were investigated by Shaikh and his coauthors (see, \cite{DBS00}, \cite{DSB03}, \cite{HMS10}, \cite{SH08}, \cite{SH08a}, \cite{SH09}, \cite{SH09a}, \cite{SH09b}, \cite{SJ07}, \cite{SJ07a}, \cite{SJ07b}, \cite{SJ07c}, \cite{SJE08}, \cite{SSH08} and also references therein). Recently in \cite{SK12} Shaikh and Kundu studied the characterization of warped product $WS_n$ and $WRS_n$.\\
\indent For a $(0,4)$-tensor $H$ and a symmetric $(0,2)$-tensor $E$ we can define two endomorphisms $X\wedge_E Y$ and $\mathscr{H}(X,Y)$ by (\cite{DDHKS00}, \cite{DGHS11}, \cite{Gl08})
\begin{eqnarray*}
(X\wedge_E Y)X_1 = E(Y,X_1)X-E(X,X_1)Y \ \ \mbox{ and } \ \ 
\mathscr{H}(X,Y)X_1 = \mathcal H(X,Y)X_1
\end{eqnarray*}
respectively, where $\mathcal H$ is the corresponding $(1,3)$-tensor of $H$. We note that the corresponding endomorphism of $R$ is given by $\mathscr R(X,Y) =[\nabla_X,\nabla_Y] - \nabla_{[X,Y]}$, where the square bracket is the Lie bracket over $\chi(M)$.\\
\indent For a $(0,k)$-tensor $T$, $k\geq 1$, a $(0,4)$-tensor $H$ and a symmetric $(0,2)$-tensor $E$ we define two
$(0,k+2)$-tensors $H \cdot T$, and $Q(E, T)$ by (\cite{DDHKS00}, \cite{DGHS11}, \cite{DGPSS11}, \cite{Gl08})
\begin{eqnarray*}
&&(H\cdot T)(X_1,\cdots, X_k,X,Y)=(\mathcal{H}(X,Y)\cdot T)(X_1,\cdots, X_k)\\
&&=-T(\mathcal{H}(X,Y)X_1,X_2,\cdots,X_k)-\cdots-T(X_1,\cdots,X_{k-1},\mathcal{H}(X,Y)X_k),\\
&&Q(E, T)(X_1,\cdots, X_k;X,Y)=((X\wedge_E Y)\cdot T)(X_1,\cdots, X_k)\\
&&=-T((X\wedge_E Y)X_1,X_2,\cdots,X_k)-\cdots-T(X_1,\cdots,X_{k-1},(X\wedge_E Y)X_k).
\end{eqnarray*}
Putting in the above formulas $T = R, S, C, P, W, K$ and $E = g$ or $S$, 
we obtain the tensors: $R \cdot R$, $R \cdot S$, $R \cdot C$, $R \cdot P$, $R \cdot W$, $R \cdot K$, $C \cdot R$, $C \cdot S$, $C \cdot C$, $C \cdot P$, $C \cdot W$, $C\cdot K$, $P \cdot R$, $P \cdot C$, $P \cdot P$, $P \cdot W$, $P \cdot K$, $Q(g, R)$, $Q(g, S)$, $Q(g, C)$, $Q(g, P)$, $Q(g, W)$, $Q(g, K)$, $Q(S, R)$, $Q(S, C)$, $Q(S, P)$, $Q(S, W)$, $Q(S, K)$ etc. The tensor $Q(E,T)$ is called the Tachibana tensor of the tensors $E$ and $T$, or the Tachibana tensor for short \cite{DGPSS11}.\\
\indent If the manifold $M$ satisfies the condition $R \cdot R = 0$ (resp., $R \cdot S = 0$, $R \cdot C = 0$, $R \cdot W = 0$, $R \cdot P = 0$, $R \cdot K = 0$), then it is called semisymmetric (\cite{Ca46}, \cite{Sz82}) (briefly, $SS_n$) (resp., Ricci semisymmetric (briefly, $RSS_n$), conformally semisymmetric (briefly, $CSS_n$), projectively semisymmetric, concircularly semisymmetric, conharmonically semisymmetric).\\ 
\indent The manifold $M$ is said to be pseudosymmetric in sense of Deszcz or simply Deszcz pseudosymmetric (resp., Ricci pseudosymmetric) \cite{DGr87} if the condition
\begin{equation}
\label{4.3.002}
R \cdot R = L_R\, Q(g,R) \ \ (\mbox{resp.,} R \cdot S = L_{S}\, Q(g,S))
\end{equation}
holds on $U_R = \{ x\in M :\left(R -\frac{\kappa}{n(n-1)}G\right)_x \neq 0 \}$ (resp., $U_S = \{ x \in M: \left(S - \frac{r}{n}\, g\right)_x \neq 0 \}$) for a function $L_R$ (resp., $L_S$) on $U_R$ (resp., $U_S$).\\
\indent The manifold $M$, $n \geq 4$, is said to be a manifold with pseudosymmetric Weyl conformal curvature tensor (\cite{Des91}, \cite{DGr98}) (resp., pseudosymmetric Weyl projective curvature tensor) if
\begin{equation}\label{4.3.012}
C \cdot C = L_{C}\, Q(g,C) \ \ (\mbox{resp., } P \cdot P = L_{P}\, Q(g,P))
\end{equation}
on $U_C = \{ x \in M : C \neq 0\ \mbox{at}\ x \}$ (resp., $U_P = \{ x \in M : Q(g,P) \neq 0\ \mbox{at}\ x \}$), where $L_{C}$ (resp., $L_P$) is a function on $U_C$ (resp., $U_P$). We note that $U_{C} \subset U_{R}$.\\
\indent Again $M$ is said to be Ricci-generalized pseudosymmetric (\cite{DD91}, \cite{DD91a}) if at every point of $M$, the tensor $R\cdot R$ and the Tachibana tensor $Q(S,R)$ are linearly dependent. Hence $M$ is Ricci-generalized pseudosymmetric if and only if 
\begin{equation}\label{4.3.00777}
R\cdot R=L\, Q(S,R)
\end{equation}
holds on  $\{x\in M : Q(S,R) \neq 0\ \mbox{at}\ x \}$, where $L$ is some function on this set. An important subclass of Ricci-generalized pseudosymmetric manifolds is formed by the manifolds realizing the condition (\cite{DD91}, \cite{DGr98})
\begin{equation}\label{4.3.00444}
R \cdot R=Q(S,R). 
\end{equation}
\indent The manifold $M$, $n \geq 4$, is said to be a manifold with Ricci generalized pseudosymmetric Weyl conformal curvature tensor (\cite{Des91}, \cite{DGr98}) (resp., Ricci generalized pseudosymmetric Weyl projective curvature tensor) if
\begin{equation}
C \cdot C = L_1\, Q(S,C) \ \ (\mbox{resp., } P \cdot P = L_2\, Q(S,P))
\end{equation}
on $U_1 = \{ x \in M : Q(S,C) \neq 0\ \mbox{at}\ x \}$ (resp., $U_2 = \{ x \in M : Q(S,P) \neq 0\ \mbox{at}\ x \}$), where $L_1$ (resp., $L_2$) is a function on $U_1$ (resp., $U_2$).\\
\indent A symmetric $(0, 2)$-tensor $E$ on $M$ is called Riemann compatible or simply $R$-compatible (\cite{MM12}, \cite{MM12a}) if 
\[
R(\mathcal E X_1, X,X_2,X_3) + R(\mathcal E X_2, X,X_3,X_1) + R(\mathcal E X_3, X,X_1,X_2) = 0
\]
holds $\forall~ X, X_1, X_2, X_3 \in \chi(M)$, where $\mathcal E$ is the endomorphism on $\chi(M)$ defined as $g(\mathcal E X_1, X_2) = E(X_1, X_2)$.\\
Again a vector field $Y$ is called $R$-compatible (also known as Riemann compatible \cite{MM13}) if
$$\Pi(X_1)R(Y,X,X_2,X_3) + \Pi(X_2)R(Y,X,X_3,X_1) + \Pi(X_3)R(Y,X,X_1,X_2) = 0$$
holds $\forall~ X, X_1, X_2, X_3 \in \chi(M)$, where $\Pi$ is the corresponding 1-form of $Y$, i.e., $\Pi(X) = g(X,Y)$. We note that if $Y$ is Riemann compatible, then $\Pi\otimes \Pi$ is a Riemann compatible tensor. Similarly for other curvature tensors we have conformal compatibility (also known as Weyl compatibility see, \cite{DGJPZ13} and \cite{MM13}), concircular compatibility, conharmonic compatibility for both a vector and a $(0,2)$-tensor.\\
\indent Recently Mantica and Suh (\cite{MS12}, \cite{MS13}, \cite{MS14} and \cite{MS12b}) obtained a necessary and sufficient condition for the recurrency of the curvature $2$-forms $\Omega_{(R)l}^m= R_{jkl}^m dx^j \wedge dx^k$ (\cite{Be87}, \cite{LR89}) and Ricci 1-forms $\Lambda_{(S)l} = S_{lm} dx^m$, where $\wedge$ indicates the exterior product. They showed that $\Omega_{(R)l}^m$ are recurrent (i.e., $\mathcal D \Omega_{(R)l}^m  = A \wedge \Omega_{(R)l}^m$, $\mathcal D$ is the exterior derivative and $A$ is the associated $1$-form) if and only if
\beb\label{man}
&&(\nabla_{X_1} R)(X_2,X_3,X,Y)+(\nabla_{X_2} R)(X_3,X_1,X,Y)+(\nabla_{X_3} R)(X_1,X_2,X,Y) =\\
&&\hspace{1in} A(X_1) R(X_2,X_3,X,Y) + A(X_2) R(X_3,X_1,X,Y)+ A(X_3) R(X_1,X_2,X,Y)
\eeb
and $\Lambda_{(S)l}$ are recurrent (i.e., $\mathcal D \Lambda_{(S)l}  = A \wedge \Lambda_{(S)l}$) if and only if
$$(\nabla_{X_1} S)(X_2,X) - (\nabla_{X_2} S)(X_1,X) = A(X_1) S(X_2,X) - A(X_2) S(X_1,X)$$
for an 1-form $A$. If on a manifold the curvature 2-forms are recurrent for $R$ and $C$, then we denote such a manifold by $K_n(2)$ and $CK_n(2)$ respectively. Again if the Ricci 1-forms are recurrent then it is denoted by $RK_n(1)$.
\section{\bf Existence of $HGK_n$}
\begin{figure}[h]
{\epsfxsize=15cm\epsfysize=8cm\epsfbox{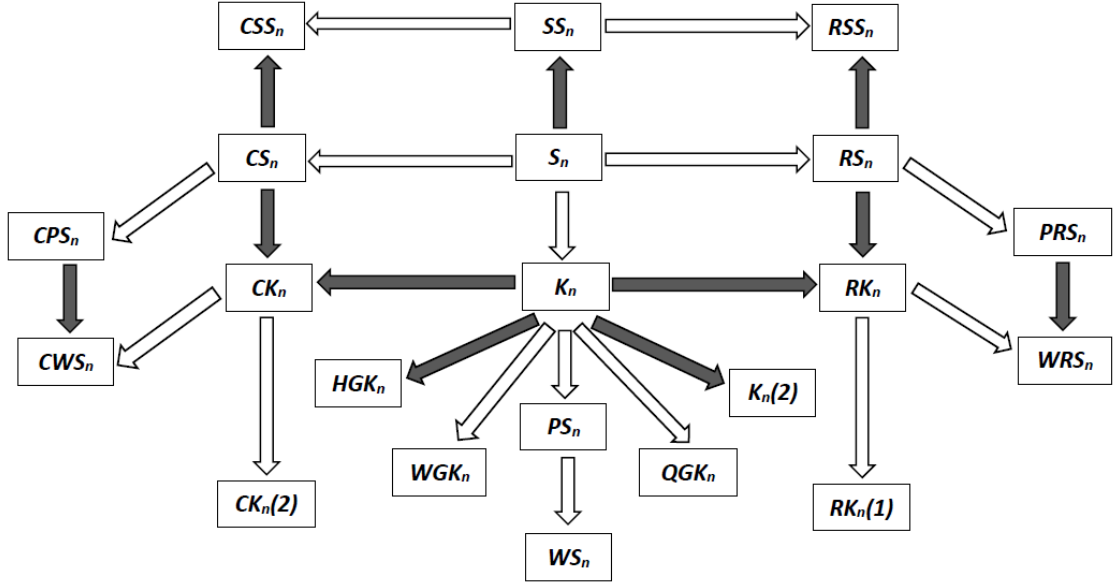}}\vspace{-0.0in}
\caption{Some generalizations of local symmetry}
\label{Fig1} 
\end{figure}

Let $M$ be an open connected subset of $\mathbb{R}^4$ such that $x^1,x^2,x^3,x^4 > 0$ endowed with the Lorentzian metric
\be\label{3.1}
ds^2 = g_{ij}dx^i dx^j = e^{x^1+x^3} (dx^1)^2 + 2 dx^1 dx^2 + (dx^3)^2 + e^{x^1} (dx^4)^2, \ \ i, j = 1, 2, 3, 4.
\ee
Then the  non-zero components of the Christoffel symbols of second kind (upto symmetry) are given by:
$$\Gamma^{2}_{11} = -\Gamma^{3}_{11} = \Gamma^{2}_{13} =\frac{e^{x^1+x^3}}{2}, \ \Gamma^{4}_{14}=\frac{1}{2}, \ \Gamma^{2}_{44}=-\frac{e^{x^1}}{2}.$$
The non-zero components of curvature tensor and Ricci tensor (upto symmetry) are given by:
\be\label{3.2}
R_{1313}=-\frac{1}{2} e^{x^1+x^3}, \ \ \ R_{1414}=-\frac{e^{x^1}}{4}, \ \ \ \ \ S_{11}=\frac{1}{4}(1+2 e^{x^1+x^3}).
\ee
The scalar curvature of this metric is given by $r = 0$. Again the non-zero components $R_{hijk,l}$ and $S_{ij,l}$ of the covariant derivatives of curvature tensor and Ricci tensor (upto symmetry) are given by:
\be\label{3.3}
R_{1313,1}=-\frac{e^{x^1+x^3}}{2}=R_{1313,3}, \ \ \ \ S_{11,1}= S_{13,3}=-\frac{e^{x^1+x^3}}{2}.
\ee
The non-zero components of the conformal curvature tensor and its covariant derivative (upto symmetry) are given by:
\be\label{3.4}\left\{\begin{array}{c}
C_{1313} = \frac{1}{8} (1 - 2 e^{x^1+x^3}), \ \ \ C_{1414} = -\frac{1}{8} e^{x^1}(1 - 2 e^{x^1+x^3}),
\vspace{0.1in}\\
C_{1313, 1} = C_{1313, 3} = -\frac{e^{x^1+x^3}}{4}, \ \ \ C_{1414, 1}= C_{1414, 3} =\frac{e^{2x^1+x^3}}{4}.\end{array}\right.
\ee
The non-zero components of the conharmonic curvature tensor and its covariant derivative (upto symmetry) are given by:
\be\label{3.5}\left\{\begin{array}{c}
K_{1313} = \frac{1}{8} (1 - 2 e^{x^1+x^3}), \ \ \ K_{1414} = -\frac{1}{8} e^{x^1}(1 - 2 e^{x^1+x^3}),
\vspace{0.1in}\\
K_{1313, 1} = K_{1313, 3} = -\frac{e^{x^1+x^3}}{4}, \ \ \ K_{1414, 1}= K_{1414, 3} =\frac{e^{2x^1+x^3}}{4}.\end{array}\right.
\ee
The non-zero components of the projective curvature tensor and its covariant derivative (upto symmetry) are given by:
\be\label{3.6}\left\{\begin{array}{c}
P_{1211}=\frac{1}{12} \left(2 e^{x^1+x^1}+1\right), \ \ P_{1313}=\frac{1}{12} \left(1-4 e^{x^1+x^1}\right), \ \ P_{1331}=\frac{1}{2} e^{x^1+x^1}\\
P_{1414}=\frac{1}{6} e^{x^1} \left(e^{x^1+x^1}-1\right), \ \ \ P_{1441}=\frac{e^{x^1}}{4}\\
\vspace{0.1in}
P_{1211,1}= P_{1211,3}= -P_{1313,1}= -P_{1313,3}= P_{1331,1}= P_{1331,3}=\frac{1}{2} e^{x^1+x^3}\\
P_{1414,1}= P_{1414,3}=\frac{1}{6} e^{2x^1+x^3}.$$
\end{array}\right.
\ee
The non-zero components of the concircular curvature tensor and its covariant derivative (upto symmetry) are given by:
\be\label{3.7}
W_{1313}=-\frac{1}{2} e^{x^1+x^3}, \ \ W_{1414}=-\frac{e^{x^1}}{4} \ \mbox{ and } \ W_{1313,1}= W_{1313,3}=-\frac{1}{2} e^{x^1+x^3}.
\ee
In view of \eqref{3.2}--\eqref{3.7}, it can be easily shown that the tensors $R\cdot R$, $C\cdot R$, $P\cdot R$, $P\cdot S$, $P\cdot C$, $P\cdot Z$, $P\cdot K$, $W\cdot R$ and $K\cdot R$ are identically zero. Thus from \cite{SK14} we can conclude that $R\cdot H$, $C\cdot H$, $W\cdot H$ and $K\cdot H$ are all identically zero for $H = C, W$ and $K$.\\
The non-zero components of $P \cdot P$ and $Q(S, P)$ are (upto symmetry) given below:
\be\label{3.8}
(P \cdot P)_{131131}= \frac{1}{144} \left(1+2 e^{x^1+x^3}\right)^2, \ \ \ (P \cdot P)_{141141}=\frac{1}{144} e^{x^1} \left(1+2 e^{x^1+x^3}\right)^2.
\ee
\be\label{3.9}
Q(S, P)_{141114} = e^{x^1} Q(S, P)_{131113} = \frac{e^{x^1}}{48}(1+2e^{x^1+x^3})^2.
\ee
\vspace{0.1in}\\
In terms of local coordinates the defining condition of a $HGK_{n}$ can be written as
\be\label{3.10}
R_{hijk,\, l}=A_{l}R_{hijk}+B_{l}[S_{hk}g_{ij} + S_{ij}g_{hk} - S_{hj}g_{ik} - g_{hj}S_{ik}], \ \ 1 \le h, i, j, k, l \le n.
\ee
Then by virtue of (\ref{3.1}), (\ref{3.2}), (\ref{3.3}) and (\ref{3.10}) it is easy to check that the manifold is a $HGK_4$ with associated 1-forms
\begin{equation}\left\{\begin{array}{ccc}\label{3.11}
&A_{i}=\left\{\begin{array}{ccc}&-\frac{2e^{x^1+x^3}}{1-2e^{x^1+x^3}}&\ \ \ \ \mbox{for}\ \ \ i=1, 3\\
																		& 0&\ \ \ \ \mbox{otherwise,}
									\end{array}\right.&\\
\\
&B_{i}=\left\{\begin{array}{ccc}
					&\frac{2e^{x^1+x^3}}{1-4e^{2x^1+2x^3}}&\ \ \ \ \mbox{for}\ \ \ i=1, 3\\
					& 0&\ \ \ \ \mbox{otherwise.}\end{array}\right.&
								\end{array}\right.
\end{equation}
\indent If we consider the 1-form $\eta =\{1, 0, 0, 0\}$ and the scalar $\beta = \frac{1}{4}(1+2 e^{x^1+x^3})$ then from (\ref{2.1}) and (\ref{3.2}) it follows that the manifold satisfies $S=\beta \eta\otimes \eta$. Thus the manifold is Ricci simple.\\
\indent Now it is easy to check that $\eta$ is parallel, i.e., $\nabla \eta = 0$, hence from Proposition \ref{pr2.1}, the manifold is Ricci recurrent. Also by virtue of \eqref{3.2} and \eqref{3.3}, we can easily check that the manifold is Ricci recurrent with the associated 1-form
\[A_{i}=\left\{\begin{array}{ccc}&\frac{2e^{x^1+x^3}}{1+2e^{x^1+x^3}}&\ \ \ \ \mbox{for}\ \ \ i=1, 3\\
& 0&\ \ \ \ \mbox{otherwise.}\end{array}\right.
\]
Now from (\ref{3.5}) it follows that the manifold is conharmonically recurrent with the associated 1-form
\[A_{i}=\left\{\begin{array}{ccc}&-\frac{2e^{x^1+x^3}}{1-2e^{x^1+x^3}}&\ \ \ \ \mbox{for}\ \ \ i=1, 3\\
& 0&\ \ \ \ \mbox{otherwise.}\end{array}\right.
\]
Now from the values of $R$ it can be easily checked that $M$ satisfies the condition
$$\Pi(X_1) R(X_2,X_3,X,Y) + \Pi(X_2) R(X_3,X_1,X,Y)+ \Pi(X_3) R(X_1,X_2,X,Y) = 0$$
for the 1-form
\[\Pi_{i}=\left\{\begin{array}{ccc}&\frac{2e^{x^1+x^3}}{1+2e^{x^1+x^3}}&\ \ \ \ \mbox{for}\ \ \ i=1\\
& 0&\ \ \ \ \mbox{otherwise.}\end{array}\right.
\]
Thus from \eqref{man}, it follows by virtue of second Bianchi identity that the curvature 2-forms $\Omega_{(R)l}^m$ are recurrent.\\
The general form of Riemann compatible tensor $E$ with local components $E_{ij}$, of $M$ can be easily be evaluated as
$$\left(
\begin{array}{cccc}
 E_{11} & E_{12} & E_{13} & E_{14} \\
 E_{21} & 0 & 0 & 0 \\
 E_{31} & 0 & E_{33} & 2 e^{x^1+x^3} E_{43} \\
 E_{41} & 0 & E_{43} & E_{44},
\end{array}
\right)$$
where $E_{ij}$'s are arbitrary.
Again the general form of Riemann compatible vector $Y$ is\\
$\left\{0, a_1, 0, a_4 e^{-x^1}\right\}$ or $\left\{{a_2, -a_2 e^{x^1 + x^3}, a_3, 0}\right\}.$\\
Similarly we can get the conformal compatible, concircular compatible and conharmonic compatible tensors respectively given in the following form:
\beb
\left(
\begin{array}{cccc}
 E_{11} & E_{12} & E_{13} & E_{14} \\
 E_{21} & 0 & 0 & 0 \\
 E_{31} & 0 & E_{33} & E_{34} \\
 E_{41} & 0 & -E_{34} & E_{44}
\end{array}
\right),  
\left(
\begin{array}{cccc}
 E_{11} & E_{12} & E_{13} & E_{14} \\
 E_{21} & 0 & 0 & 0 \\
 E_{31} & 0 & E_{33} & 2 e^{x^1+x^3} E_{43} \\
 E_{41} & 0 & E_{43} & E_{44}
\end{array}
\right) \mbox{and} 
\left(
\begin{array}{cccc}
 E_{11} & E_{12} & E_{13} & E_{14} \\
 E_{21} & 0 & 0 & 0 \\
 E_{31} & 0 & E_{33} & E_{34} \\
 E_{41} & 0 & -E_{34} & E_{44}
\end{array}
\right).
\eeb
Again the general form of conformal (resp., concircular and conharmonic) compatible vector $Y$ is $\left\{0, a_1, 0, a_4 e^{-x^1}\right\}$ or $\left\{{a_2, -a_2 e^{x^1 + x^3}, a_3, 0}\right\}.$
%
\begin{rem}
From Lemma 5.6 of \cite{SK14} it follows that if a semi-Riemannian manifold $M$ is Ricci recurrent and also conformally recurrent with same 1-form of recurrence, then $M$ is recurrent. We note that the above manifold is Ricci recurrent and also conformally recurrent but of different 1-forms of recurrence and it is not recurrent.
\end{rem}
Now we can state the following:
\begin{pr}
Let $(M^4, g)$ be a semi-Riemannian manifold equipped with the metric (\ref{3.1}). Then $(M^4, g)$ is:\\
(i) Ricci simple and hence a special quasi Einstein manifold,\\
(ii) manifold of constant scalar curvature,\\
(iii) Ricci recurrent and hence weakly Ricci symmetric and the Ricci 1-forms are recurrent $($see \cite{SDHJK15} and \cite{MS12b}$)$,\\
(iv) conharmonically recurrent and hence conformally recurrent $($see Corollary 6.2 of \cite{SK14}$)$ as well as conharmonically and conformally weakly symmetric, and the curvature 2-forms for $C$ and $K$ are recurrent,\\
(v) hyper-generalized recurrent and hence generalized Ricci recurrent,\\
(vi) semisymmetric and hence fulfills the semisymmetry conditions for $S$, $C$, $P$, $W$ and $K$,\\
(vii) weakly Ricci symmetric and weakly conharmonic symmetric such that all the  associated 1-forms are different $($This result proves the non uniqueness of associated 1-forms of weak symmetry structures, see Section 3 of \cite{SDHJK15}$)$,\\
(viii) recurrent curvature 2-form, and\\
(ix) fulfills the condition $P \cdot P = -\frac{1}{3}Q(S, P)$, i.e, Ricci generalized pseudosymmetric Weyl projective curvature tensor.
\end{pr}
From the above values of various tensors, we can state the following:
\begin{pr}\label{pr3.2}
The 4-dimensional semi-Riemannian manifold $M$ equipped with the metric given in (\ref{3.1}) is not any one the following:\\
(i) Einstein, (ii) class $\mathcal A$ or $\mathcal B$, (iii) pseudosymmetric in sense of Chaki for $R$, $C$, $P$, $W$ and $K$ (iv) pseudo Ricci symmetric, (v) quasi generalized recurrent, (vi) weakly generalized recurrent, (vii) $P\cdot P =0$, (viii) Codazzi type Ricci tensor, (ix) cyclic Ricci parallel, (x) weakly symmetric for $R$ and $Z$.
\end{pr}
It is well known that recurrent structure implies Ricci recurrent, conformally recurrent and also conharmonically recurrent but the converse is not true in each case. Again Chaki pseudosymmetry for $R$, $S$, $C$, $P$, $W$ and $K$ implies weak symmetry for $R$, $S$, $C$, $P$, $W$ and $K$ respectively but the converse cases are not true. Now from above two propositions we can state the following:
\begin{theo}\label{th3.1}
Let $(M^4, g)$ be a semi-Riemannian manifold equipped with the metric (\ref{3.1}). Then $(M^4, g)$ is\\
(i) not locally symmetric but semisymmetric,\\
(ii) not weakly symmetric but weakly Ricci symmetric as well as conformally weakly symmetric,\\ 
(iii) not recurrent but Ricci recurrent as well as conformally recurrent,\\
(iv) not recurrent but hyper-generalized recurrent,\\
(v) neither class $\mathcal{A}$ nor $\mathcal{B}$ but of constant scalar curvature,\\
(vi) not Einstein but Ricci simple i.e., special quasi-Einstein,\\
(vii) not recurrent but the curvature 2-forms are recurrent.\\
(viii) does not satisfy $P\cdot P =0$ but fulfills the condition $P\cdot P = L Q(S, P)$.
\end{theo}
\begin{rem}\label{rem3.1}
From above theorem we get an affirmative answer of \textbf{Q. 1} and we can state the following:\\
(i) semisymmetry is a proper generalization of local symmetry and the result is also true for the case $S$, $C$, $P$, $W$ and $K$,\\
(ii) Ricci recurrency and conharmonic recurrency are both proper generalizations of recurrency,\\
(iii) weak Ricci symmetry is a proper generalization of pseudo Ricci symmetry and the result is also true for conharmonic and concircular curvature tensors,\\
(iv) the classes $\mathcal{A}$ and $\mathcal{B}$ properly lie between Ricci symmetry and manifold of constant scalar curvature,\\
(v) hyper-generalized recurrency is a proper generalization of recurrency,\\
(vi) recurrency of curvature 2-forms is a proper generalization of recurrency.
\end{rem}
The above remark can be stated in a diagram as shown in Figure \ref{Fig1}, where the boxes denote the corresponding class of manifolds and the `implication' signs denote that the later one is a generalization of the previous. The above example ensures the properness of some of them, which are indicated by the dark filled ones.
\vspace{-0.0in}
\\
%
\indent Again, if we change the signature of the metric (\ref{3.1}) as
\[
ds^2 = g_{ij}dx^i dx^j = -e^{x^1+x^3} (dx^1)^2 + 2 dx^1 dx^2 + (dx^3)^2 + e^{x^1} (dx^4)^2,
\]
\[
ds^2 = g_{ij}dx^i dx^j = e^{x^1+x^3} (dx^1)^2 + 2 dx^1 dx^2 - (dx^3)^2 - e^{x^1} (dx^4)^2,
\]
\[
ds^2 = g_{ij}dx^i dx^j = -e^{x^1+x^3} (dx^1)^2 + 2 dx^1 dx^2 - (dx^3)^2 - e^{x^1} (dx^4)^2,
\]
\[
ds^2 = g_{ij}dx^i dx^j = e^{x^1+x^3} (dx^1)^2 + 2 dx^1 dx^2 - (dx^3)^2 + e^{x^1} (dx^4)^2,
\]
\[
ds^2 = g_{ij}dx^i dx^j = -e^{x^1+x^3} (dx^1)^2 + 2 dx^1 dx^2 - (dx^3)^2 + e^{x^1} (dx^4)^2,
\]
\[
ds^2 = g_{ij}dx^i dx^j = -e^{x^1+x^3} (dx^1)^2 + 2 dx^1 dx^2 + (dx^3)^2 - e^{x^1} (dx^4)^2,
\]
\[
ds^2 = g_{ij}dx^i dx^j = e^{x^1+x^3} (dx^1)^2 + 2 dx^1 dx^2 + (dx^3)^2 - e^{x^1} (dx^4)^2,
\]
$i, j = 1, 2, 3, 4$, then it can be easily shown that the geometric properties presented in Theorem \ref{th3.1} remain unchanged for the above metrics. We note that the first three of the above metrics are Lorentzian and remaining are semi-Riemannian.
\begin{rem}
Thus we get the existence of a new class of semi-Riemannian manifolds having the geometric properties presented in Theorem \ref{th3.1}. It is similar to show that the following metrics (also with different signatures) belong to this class:
\be\label{3.12}
ds^2 = e^{x^1+x^3} (dx^1)^2 + 2 dx^1 dx^2 + (dx^3)^2 + f(x^1) (dx^4)^2,
\ee
\be\label{3.13}
ds^2 = x^1 x^3 (dx^1)^2 + x^1 (dx^2)^2 + 2 dx^1 dx^3 + f(x^1) (dx^4)^2,
\ee
\be\label{3.14}
ds^2 = e^{x^1+x^3} (dx^1)^2 + 2 dx^1 dx^2 + (dx^3)^2 + f(x^1) (dx^4)^2 + f(x^1) \delta_{ab}dx^a dx^b,
\ee
\be\label{3.15}
ds^2 = x^1 x^3 (dx^1)^2 + x^1 (dx^2)^2 + 2 dx^1 dx^3 + f(x^1) (dx^4)^2 + f(x^1) \delta_{ab}dx^a dx^b
\ee
where $f$ is a smooth function and $\delta_{ab}$ denotes the Kronecker delta, $5 \le a, b \le n$.
\end{rem}
\noindent\textbf{Conclusion:} From the above results and discussion we get the existence of a new class of semi-Riemannian manifolds which are hypergeneralized recurrent, Ricci recurrent, conharmonically recurrent, conformally recurrent, manifold of recurrent curvature 2-forms, Ricci simple (thus special quasi Einstein), weakly Ricci symmetric, conformally weakly symmetric, conharmonically weakly symmetric, semisymmetric, Ricci semisymmetry and conformally semisymmetric and fulfills the pseudosymmetric type condition $P \cdot P = -\frac{1}{3}Q(S, P)$ but do not satisfy any one of the conditions (i)-(x) of Proposition \ref{pr3.2}. Hence the existence of a generalized class of recurrent manifolds, viz., hyper generalized recurrent manifolds is ensured by the metric.\\
\indent\newline
\noindent\textbf{Acknowledgement.} We have made all the calculations by Wolfram Mathematica. The last named author gratefully acknowledges to CSIR, New Delhi, India (File No. 09/025 (0194)/2010-EMR-I) for the financial support.

\end{document}